\documentclass[12pt]{article}
\usepackage{amsmath}
\usepackage{amscd}
\usepackage{amsthm}
\usepackage{amssymb}
\usepackage[latin1]{inputenc}
\usepackage[all]{xy}
\def\email#1{{\tt #1}}
\def\SetSize{\fontsize{12}{14.4}\selectfont}

\newtheoremstyle{theorem}
  {}
  {}
  {\itshape}
  {}
  {}
  {.}
  {.5em}
  {}

\newtheoremstyle{definition}
  {}
  {}
  {}
  {}
  {\MakeUppercase}
  {.}
  {.5em}
  {}

\theoremstyle{theorem}
\newtheorem{theorem}{THEOREM}[section]

\newtheorem{corollary}[theorem]{COROLLARY}
\newtheorem{proposition}[theorem]{PROPOSITION}

\theoremstyle{definition}
\newtheorem{remark}[theorem]{REMARK}
\newtheorem{remarks}[theorem]{REMARKS}

\newtheorem{exercise*}[enumiv]{EXERCISE}

\let\phi=\varphi

\def\Cok{\operatorname{Cok}}

\def\Rad{\operatorname{Rad}}

\def\grade{\operatorname{grade}}
\def\Hom{\operatorname{Hom}}

\def\rank{\operatorname{rank}}

\def\Ker{\operatorname{Ker}}
\def\Coker{\operatorname{Coker}}

\def\Im{\operatorname{Im}}

\let\oldbigwedge\bigwedge
\def\BIGwedge{{\textstyle\oldbigwedge}}
\def\medwedge{{\scriptstyle\oldbigwedge}}
\def\bigwedge{\mathchoice{\BIGwedge}{\BIGwedge}{\medwedge}{}}

\let\iso=\cong

\let\epsilon=\varepsilon
\let\tilde=\widetilde

\begin{document}
\thispagestyle{empty} \vspace*{1.5in}
{\fontsize{14}{16.8}\selectfont \noindent \uppercase{Koszul
Bicomplexes and generalized Koszul complexes in projective
dimension one}\par}

\SetSize \vspace{2\baselineskip}

\noindent\uppercase{Bogdan Ichim}, Universit\"at Osnabr\"uck, FB
Mathematik/Informatik, 49069 Osna\-br\"uck, Germany,
\email{bogdan.ichim@mathematik.uni-osnabrueck.de}\\
Institute of Mathematics, C.P. 1-764, 70700 Bucharest, Romania,\\
\email{bogdan.ichim@imar.ro}
\\[1\baselineskip]
\noindent\uppercase{Udo Vetter}, Universit\"at Oldenburg, Institut
für Mathematik, 26111 Oldenburg, Germany,
\email{vetter@\allowbreak mathematik.uni-oldenburg.de}

\vspace{4\baselineskip plus 1 \baselineskip minus 1\baselineskip}

\begin{abstract}
We describe Koszul type complexes associated with a linear map
from any module to a free module, and vice versa with a linear map
from a free module to an arbitrary module, generalizing the
classical Koszul complexes. Given a short complex of finite free
modules, we assemble these complexes to what we call Koszul
bicomplexes. They are used in order to investigate the homology of
the Koszul complexes in projective dimension one. As in the case
of the classical Koszul complexes this homology turns out to be
grade sensitive. In a special setup we obtain necessary conditions
for a map of free modules to be lengthened to a short complex of
free modules.
\end{abstract}

\section{Introduction}
The paper deals with the following problem: let $\mathcal
F\overset \chi\to\mathcal G\overset\lambda\to \mathcal
 H$ be a complex of finite free modules over a
noetherian ring $R$; in which way does $\grade I_\lambda$ depend
on
 $\grade I_\chi$ and on the ranks of $\mathcal{F}$,
$\mathcal{G}$, $\mathcal{H}$? (By $I_\lambda$ we denote the ideal
of maximal minors of $\lambda$.) If, for example, $\rank
\mathcal{F}=1$, $\rank \mathcal{G}=n$, and $\chi$ is given by a
regular sequence $x_1,\ldots,x_n$ in $R$, i.e.
$\chi(1)=(x_1,\ldots,x_n)$, then it was proved in [BV4] that
$\grade I_\lambda=n$  is possible if and only if $\rank
\mathcal{H}=1$ and $n$ is even.

Being far away from a satisfactory solution to the the general
problem, we shall report here on some approach.

It turns out that there is a close connection to the behavior of
the homology of the generalized Koszul complex associated with the
induced map $\bar\lambda:\Cok\chi\to\mathcal H$. Here we must
clear up what we mean by a generalized Koszul complex since even
the terminology concerning the {\it classical Koszul complex} is
not standardized. We refer to the accesses given in [BO2], § 10 or
in [BH], 1.6; in [E], Ch. 17 there is a comprehensive presentation
of the `dual' approach.

The Eagon-Northcott family of complexes $\mathcal{C}^{t}(\psi)$
(see [E], A2.6 or [BV2], 2.C), associated with a linear map
$\psi:G\to F$ of finite free modules $G$ and $F$ such that $\rank
G\ge \rank F$, generalizes the Koszul complex in a special case.
The homology of $\mathcal{C}^{t}(\psi)$ is well-understood. As in
the classical case it is grade sensitive with respect to the ideal
$I_\psi$.

More generally we consider linear maps $\psi:G\to F$, where only
$F$ has to be finite and free (weaker assumptions on $F$ are
possible). We construct a family of complexes
$\mathcal{C}_\psi(t)$ associated with $\psi$ which generalizes
both the Eagon-Northcott family of complexes, and the classical
Koszul complex. There is a similar `dual' construction of a family
of complexes $\mathcal{D}_\phi(t)$ for a map $\phi:H\to G$ from a
finite free module $H$ to an arbitrary module $G$. The complexes
just mentioned are the {\it generalized Koszul complexes}.

If $H\stackrel{\phi}{\to} G\stackrel{\psi}{\to} F$ is a complex,
we can compose $\mathcal{C}_\psi(t)$ and $\mathcal{D}_\phi(t)$ to
the Koszul bicomplex $\mathcal{K}_{.,.}(t)$. Setting $H=\mathcal
H^*,G=\mathcal G^*,F=\mathcal F^*$ and
$\phi=\lambda^*,\psi=\chi^*$, we use it to investigate the
homology of $\mathcal{C}_{\bar\lambda}(t)$. Assume that
$\rank\mathcal F\le \rank\mathcal G$. Then the most satisfactory
result (see Theorem \ref{MaximalCaseHomology}) is obtained if
$\grade I_\chi$ has the greatest possible value $\rank\mathcal
G-\rank\mathcal F+1$. It covers a result of Migliore, Nagel and
Peterson (see Proposition 5.1 in [MPN]) who proved it partially
for Gorenstein rings $R$, using local cohomology. It also
generalizes Theorem 5 in [BV4]. In the case under consideration a
full answer to our initial problem is contained in Theorem
\ref{EXTENSION H-B}, a generalization of Corollary 3 in [BV1].

What can be deduced if $\grade I_\chi$ has not the greatest
possible value (but is not too small)? In this case Theorem
\ref{EXTENSION H-B submaximal} provides some answer to the
question in the beginning. As a consequence we derive a purely
numerical criterion for the non-vanishing of a product of matrices
(Corollary \ref{MatrixCriterion}).

The paper is based on results of Bruns and Vetter (see [BV3] and
[BV4]). They study the homology of the Koszul complex associated
with a linear form on a module of projective dimension one, using
a Koszul bicomplex construction obtained from a Koszul complex and
certain Eagon-Northcott complexes. The idea to build and link
Koszul bicomplexes appears also in the paper [HM] of Herzog and
Martsinkovsky (see in particular the gluing construction for the
residue field of a complete intersection).

The article for the most part consists of considerations and
results contained in the thesis [I] of the first author.

\section{Preliminaries}\label{Preliminaries}

In the paper $R$ always denotes a commutative ring. We adopt the
following convention for the {\it graded dual} of a graded module:
if $M=\oplus_{i\ge 0}M_i$ is a graded $R$-module, we shall
improperly write $M^*$ for the graded dual of $M$, that is $$
M^*=M^*_{gr}=\oplus_{i\ge 0}(M_i)^*. $$ We use it mainly in the
case in which $M$ is the symmetric algebra $S(N)$ or the exterior
algebra $\bigwedge N$ of an $R$-module $N$. Then
$$
S(N)^*=\oplus_{i\ge 0}S_i(N)^*,\quad (\bigwedge N)^*=\oplus_{i\ge
0}(\bigwedge^iN)^*.
$$
Correspondingly $M^{**}$ means the graded bidual
$(M^*_{gr})^*_{gr}$ of $M$.
\smallskip

Let $G$ be an $R$-module. The natural graded algebra homomorphism
$$
\theta:\bigwedge G^*\to (\bigwedge G)^*,
$$
is given by
$$
\theta(y^*_1\wedge\ldots\wedge y^*_p)(y_1\wedge\ldots\wedge
y_p)=\det (y^*_j(y_i))
$$
for all $y_1,\ldots,y_p\in G$ and $y^*_1,\ldots,y^*_p\in G^*$. If
$G$ is finitely generated and projective, then $\theta$ is an
isomorphism.
\smallskip

Beside the natural left $\bigwedge G$-module structure, the
exterior algebra $\bigwedge G$ carries a natural structure as a
right $\bigwedge G^*$-module: we set
$$
y_1\wedge\ldots\wedge y_n\leftharpoonup y^*_1\wedge\ldots\wedge
y^*_p=\sum_{\sigma}\epsilon(\sigma) \det_{1\le i,j\le p}
(y^*_j(y_{\sigma(i)}))y_{\sigma(p+1)}\wedge\ldots\wedge
y_{\sigma(n)}
$$
for $y_1,\ldots ,y_n\in G$ and $y^*_1,\ldots ,y^*_p\in G^*$, where
$\sigma$ runs through the set of permutations of $n$ elements
which are increasing on the intervals $[1,p]$ and $[p+1,n]$. In
fact, the operation $\leftharpoonup$ makes $\bigwedge G$ a right
$\bigwedge G^*$-module (see [BO1], Ch. III, § 11.9 for details).

An easy calculation shows that an element  $z^*\in G^*$ acts like
an antiderivation on $\bigwedge G$ in the sense that
\begin{equation}\label{antiderivation}
(x\wedge y)\leftharpoonup z^*=(x\leftharpoonup z^*)\wedge
y+(-1)^{deg x}x\wedge(y\leftharpoonup z^*)
\end{equation}
for elements $x,y\in \bigwedge G$, $x$ homogeneous. By induction
we obtain immediately the following {\it associativity formula}:

\begin{proposition}\label{DIAGRAMCOMMUTATIVITY} Let $x_i\in G$ for $i=1,\ldots,
l$, and let $z^*_j\in G^*$ for $j=1,\ldots,p$. Assume that
$z^*_j(x_i)=0$ for all $i$, $j$. Then
$$
x_1\ldots x_l\wedge( y_1\ldots y_n\leftharpoonup z^*_1\ldots
z^*_p)=(-1)^{lp}(x_1\ldots x_l\wedge y_1\ldots y_n)\leftharpoonup
z^*_1\ldots z^*_p
$$
for all $y\in \bigwedge G$.
\end{proposition}

Furthermore we shall make use of
\begin{proposition}\label{DUALIZATION_COMMUTATIVITY}Let $y^*\in \bigwedge G^*$ be homogeneous of degree $p$.
Let $y^*\wedge$ denote the left multiplication by $y^*$ on
$\bigwedge G^*$, and $\leftharpoonup y^*$ the right multiplication
by $y^*$ on $\bigwedge G$. Then the diagram
$$
\CD
\bigwedge^q G^*@>{y^*\wedge}>>\bigwedge^{p+q} G^*\\
@V\theta VV @V\theta VV\\
(\bigwedge^q G)^*@>{(\leftharpoonup y^*)^*}>>(\bigwedge^{p+q}
G)^*,
\endCD
$$
is commutative.
\end{proposition}

\begin{proof}Let $x\in \bigwedge^{p+q}G$ and $z^*\in\bigwedge^q G^*$. Then
$$
\theta(y^*\wedge z^*)(x)=x\leftharpoonup(y^*\wedge
z^*)=(x\leftharpoonup y^*)\leftharpoonup z^*
=\theta(z^*)(x\leftharpoonup y^*)=(\leftharpoonup
y^*)^*(\theta(z^*))(x).
$$
\end{proof}
The last paragraph of this preparing section is devoted to the
divided power algebra $D(H)$ of a free $R$-module $H$ of finite
rank. For the definition we refer to [E], A2.4. Here we only
explain the natural structure of $D(H)$ as a graded
$S(H^*)$-module which we shall use in the next section.

The natural structure of $S(H)^*$ as a graded $S(H)$-module is
defined by
$$(ay^*)(x)=y^*(ax)$$
for all $a,x\in S(H)$ and all $y^*\in S(H)^*$. Replace $H$ by
$H^*$, use the canonical equality $H=H^{**}$, and set $$(a^*
c)(b^*)=\alpha(c)(a^*b^*)$$ for all $a^*,b^*\in S(H^*)$ and $c\in
D(H)$ where the graded algebra isomorphism $\alpha:D(H)\to
S(H^*)^*$ is given by the formula
$$\alpha(h^{(k)})(\prod_i (h^*_i)^{k_i})=
\begin{cases}
\qquad 0 & \text{if $\quad\sum k_i\neq k$}\\ \prod\limits_{i}
(h^*_i(h))^{k_i} & \text{if $\quad\sum k_i=k$},
\end{cases}
$$
$h\in H$, $h^*_i\in H^*$ (see [BE2], 1 A')). The interaction
between the $S(H^*)$-module structure of $D(H)$ and the algebra
structure of $D(H)$ is described by the following result which we
quote without its easy proof.
\begin{proposition}\label{DUALACTION}If $h^*\in H^*$, then $h^*$ acts like a derivation on $D(H)$
in the sense that
\begin{enumerate}
\item[\rm{(a)}]$h^* h^{(k)}=h^{*}(h)\, h^{(k-1)}$ for all $h\in
H$, $k\ge 1$, \item[\rm{(b)}]$h^*(cd)=(h^* c)\, d+c\, (h^* d)$ for
all $c,d\in D(H)$.
\end{enumerate}
\end{proposition}
\bigskip

\section{Generalized Koszul Complexes}\label{Generalized Koszul
Complexes}

 Let
$
 H\overset\phi\to G$ and $G\overset\psi\to F
$ be homomorphisms of $R$-modules. Most of the results of this
section are true if the canonical maps
$$
\zeta_{H,G}:H^*\otimes G \to \Hom(H,G) \quad \text{and} \quad
\zeta_{G,F}:G^*\otimes F \to \Hom(G,F)
$$
are isomorphisms. This holds, in particular, if $H$ and $F$ are
finitely generated free modules. For simplicity we restrict the
presentation to this case.

We start on associating a Koszul complex with the map $\psi$. By
virtue of $\zeta_{G,F}$ we may consider $\psi$ an element of
$\bigwedge_{S(F)} (G^*\otimes S(F))$ and, using the canonical
isomorphism of $S(F)$-algebras
$$
\bigwedge_{S(F)} (G^*\otimes S(F))\iso \bigwedge G^*\otimes S(F),
$$
an element of $\bigwedge G^*\otimes S(F)$.

Let $N$ be a right $\bigwedge G^*$-module and $M$ an
$S(F)$-module. Then $N\otimes_R M$ is a right $\bigwedge
G^*\otimes S(F)$-module by $$ (n\otimes
m)\leftharpoonup(y^*\otimes f)=(n\leftharpoonup y^*)\otimes fm,$$
$n\in N,\ m\in M,\ y^*\in \bigwedge G^*,\ f\in S(F)$. If
$M=\oplus_{i}M_i$ is a graded $S(F)$-module, then the right
multiplication by $\psi$, denoted by $\partial_\psi$, maps
$N\otimes_R M_{i-1}$ to $N\otimes_R M_i$. Since $\psi^2=0$ in
$\bigwedge G^*\otimes S(F)$, we obtain a complex of $R$-modules
$$
\cdots\to N \otimes M_{i-1}\stackrel{\partial_\psi}{\to}N \otimes
M_{i}\stackrel{\partial_\psi}{\to}N \otimes M_{i+1}\to\cdots.
$$
If in particular $N=\bigwedge G$ equipped with the $\bigwedge
G^*$-module structure defined in section \ref{Preliminaries}, then
an easy calculation shows that
\begin{align}
\partial_\psi(y_1\wedge\ldots\wedge y_p\otimes m)=\sum_{j=1}^p
(-1)^{j+1}y_1\wedge\ldots\widehat{y_j}\ldots\wedge y_p\otimes
 \psi(y_j)m
\end{align}
for all $y_1,\ldots, y_p\in G$ and $m\in M$.

We are mainly interested in two special cases for $M$, namely
$M=S(F)$ where $M_i=S_i(F)$, and $M=S(F)^*$ where
$M_i=S_{-i}(F)^*$.

Let  $f_1,\ldots,f_m$ be a basis of $F$ and $f_1^*,\ldots,f_m^*$
the dual basis of $F^*$. Then $\psi$, as an element of $\bigwedge
G^*\otimes S(F)$, has the presentation $$
\psi=\sum_{j=1}^m\psi^*(f_j^*)\otimes f_j.
$$ Set $x^*=\psi^*(f^*_1)\wedge \ldots
\wedge \psi^*(f^*_m)\in \bigwedge G^*$. We obtain immediately that
$(x^*\otimes 1_{S(F)})\psi=\psi(x^*\otimes 1_{S(F)})=0$ in
$\bigwedge G^*\otimes S(F)$. Therefore we may use the right
multiplication by $x^*$
 on $\bigwedge G$, which we denote by $\nu_{\psi}$, in order to splice
the complexes $(\bigwedge G\otimes S(F)^*,\partial_\psi)$ and
$(\bigwedge G\otimes S(F),\partial_\psi)$ to a complex
$$
\cdots\to\bigwedge G \otimes
S_{1}(F)^*\stackrel{\partial_\psi}{\to}\bigwedge G\otimes
S_{0}(F)^*\stackrel{\nu_{\psi}}{\to}\bigwedge G\otimes
S_{0}(F)\stackrel{\partial_\psi}{\to}\bigwedge G \otimes
S_{1}(F)\to\cdots.
$$
This complex, as a complex of $R$-modules, splits into direct
summands
\begin{align*}
\mathcal C_\psi(t):\quad \cdots\to\bigwedge^{t+m+p} G \otimes
S_{p}(F)^*\stackrel{\partial_\psi}{\to}\cdots
\stackrel{\partial_\psi}{\to}\bigwedge^{t+m} G & \otimes
S_{0}(F)^* \stackrel{\nu_{\psi}}{\to}\bigwedge^{t} G \otimes
S_{0}(F)\stackrel{\partial_\psi}{\to}\\ &\cdots
\stackrel{\partial_\psi}{\to} \bigwedge^{0} G \otimes S_{t}(F) \to
0.
\end{align*}

\begin{remarks}\label{USUALKOSZULINDENTIFICATION} {\rm(a)} Clearly the
connection map $\nu_\psi$ in $\mathcal C_\psi(t)$ depends on the
basis we took for $F$. To avoid notational complications we ignore
this fact, the more so as it doesn't affect the homology.

\noindent{\rm(b)} We specialize to the case in which $F=R$ and
$x^*=\psi^*(1)=\psi$. One can easily check that, for all $t\in
\mathbb{Z}$, $\mathcal C_\psi(t)$ is the classical Koszul complex
associated with $\psi$.
\end{remarks}

 There are well-known generalizations of the classical
 Koszul complex in case $G$ is a finitely
generated free $R$-module, due to
 Eagon and Northcott, Buchsbaum and Rim, and others.
Suppose $G$ to be free of finite rank, and denote by
$\mathcal{C}^t(\psi)$ the complexes as introduced in [E], A2.6.1
or in [BV2], 2.C. They are constructed as follows. Set $\rank
G=n,$ and consider the direct ($R$-)summands
$$
0\to\bigwedge^{n-m-t} G \otimes S_0(F) \to \cdots
\stackrel{\partial_\psi}{\to}\bigwedge^{0} G \otimes
S_{n-m-t}(F)\to 0$$ and
$$
0\to\bigwedge^{t} G \otimes S_0(F)\to \cdots
\stackrel{\partial_\psi}{\to}\bigwedge^{0} G \otimes S_t(F)\to 0$$
of $(\bigwedge G\otimes S(F),\partial_\psi)$. We shall connect the
$R$-dual of the first complex with the second one. Let
$\omega:\bigwedge^n G\to R$ be an orientation of $G$. Then there
is a basis $y_1,\ldots,y_n$ of $G$ such that
$\omega(y_1\wedge\ldots\wedge y_n)=1$. Denote by
$\tilde\omega:\bigwedge^n G^*\to R$ the orientation of $G^*$
induced by $\omega$, i.e. $\tilde\omega(y_1^*\wedge\ldots\wedge
y^*_n)=1$ where $y_1^*,\ldots, y^*_n$ is the basis of $G^*$ dual
to $y_1,\ldots, y_n$. Let $x^*=\psi^*(f^*_1)\wedge \ldots \wedge
\psi^*(f^*_m)$ as above, and consider the map
$$\tilde {\nu}_\psi:\bigwedge^{n-m-t} G^*\to (\bigwedge^t
G^*)^*,\quad \tilde{\nu}_\psi(z^*)(y^*)=\tilde\omega(z^*\wedge y^*
\wedge x^*)$$ for all $z^*\in \bigwedge^{n-m-t} G^*,\ y^*\in
\bigwedge^{t} G^*$. Furthermore let $\mathcal{C}^t(\psi)$ be the
$R$-sequence
\begin{align*} 0\to(\bigwedge^{0} G \otimes
S_{r-t}(F))^*\stackrel{\partial^*_\psi}{\to}\cdots
\stackrel{\partial^*_\psi}{\to}&(\bigwedge^{n-m-t} G \otimes
S_{0}(F))^* \\ & \stackrel{\tilde{\nu}_{\psi}}{\to}\bigwedge^{t} G
\otimes S_{0}(F) \stackrel{\partial_\psi}{\to}\cdots
\stackrel{\partial_\psi}{\to} \bigwedge^{0} G \otimes S_{t}(F) \to
0.
\end{align*}
Here we set $(\bigwedge^{p}G)^*=\bigwedge^{p}G^*$ and
$(\bigwedge^t G^*)^*=\bigwedge^t G$ via the corresponding
canonical maps. The following proposition shows that
$\mathcal{C}^t(\psi)$ is a complex, and that there is a complex
isomorphism $ \mathcal C_\psi(t)\iso\mathcal{C}^t(\psi). $

\begin{proposition}\label{IDENTIFICATION_1}With the assumptions and the notation from above
let $$ \omega_{p}:\bigwedge^{p} G\to(\bigwedge^{n-p} G)^*,\quad
p=0,\ldots,n,$$ be the isomorphisms given by
$$
(\omega_{p}(x))(y)=\omega(x\wedge y)$$ for $ x\in\bigwedge^{p} G,\
y\in\bigwedge^{n-p} G$. Set $S_p=S_p(F),\ S_p^*=S_p(F)^*,\
\bigwedge^p=\bigwedge^p G$. Then the diagram
$$ \xymatrix{ 0\ar[r]&\bigwedge^n\otimes S_{n-m-t}^*\ar[dd]^{\pm
\omega_n\otimes 1} \ar[r]^{\qquad\partial_\psi} &
\cdots\ar[r]^{\partial_\psi}
&\bigwedge^{m+t}\ar[dd]^{\pm\omega_{m+t}}
\ar[dr]^{\nu_\psi}\\
&&&& \bigwedge^t\ar[r]^{\partial_\psi}&\cdots \ar[r]^{\partial_\psi}&S_t\ar[r]&0\\
0\ar[r]&(\bigwedge^0\otimes S_{n-m-t})^*\ar[r]^{\qquad\quad
\partial^*_\psi} &\cdots
\ar[r]^{\negthickspace\negthickspace\negthickspace
\negthickspace\negthickspace\negthickspace\partial^*_\psi}&(\bigwedge^{n-m-t}
)^* \ar[ru]_{\tilde\nu_\psi}}
$$
is commutative. In particular there are $($non-canonical$)$
complex isomorphisms
$$
\mathcal C_\psi(t)\iso\mathcal{C}^t(\psi).
$$
\end{proposition}

\begin{proof}
We start with showing that the diagram
$$
\CD
\bigwedge^{p} G\otimes  S(F)^*@>\partial_{\psi}>>\bigwedge^{p-1} G\otimes  S(F)^*\\
@V\omega_p\otimes 1 VV @VV\omega_{p-1}\otimes 1 V\\
(\bigwedge^{n-p} G)^*\otimes
S(F)^*@>\partial^*_{\psi}>>(\bigwedge^{n-(p-1)} G)^*\otimes
S(F)^*\\
\endCD
$$
is commutative or anticommutative. Let $x\in \bigwedge ^p G,\ y\in
\bigwedge^{n-(p-1)}G$, $f\in S(F)$ and $f^*\in S(F)^*$. Using
formula (2) from above, one easily shows that
\begin{equation*}
(\omega_{p-1}\otimes 1)\partial_{\psi}  (x\otimes f^*)(y\otimes f)
=(\omega\otimes f^*)\big(\partial_\psi(x\otimes f)(y\otimes
1)\big),
\end{equation*}
and
\begin{equation*}
\partial^*_{\psi}(\omega_p\otimes 1)  (x\otimes f^*)(y\otimes
f) =(\omega\otimes f^*)\big((x\otimes 1)\partial_\psi(y\otimes
f)\big).
\end{equation*}
Since
\begin{equation*}
0=\partial_\psi(x\wedge y\otimes f)=\partial_\psi(x\otimes f) (
y\otimes 1)+(-1)^p(x\otimes 1)\partial_\psi(y\otimes f),
\end{equation*}
by formula (\ref{antiderivation}) in section \ref{Preliminaries},
we are done.

It remains to prove that the `triangle' is commutative or
anticommutative, that is $\nu_\psi=\pm
\tilde\nu_\psi\omega_{m+t}$. We consider the diagram
$$
\xymatrix{\bigwedge^{m+t} G\ar[rr]^{\nu_\psi=\leftharpoonup
x^*}\ar[d]_{\omega_{m+t}}&& \bigwedge^t
G \ar[d]^{\omega_t}\\
(\bigwedge^{n-m-t}G)^*\ar[rr]^{(\leftharpoonup x^*)^*}&&
(\bigwedge^{n-t} G)^*
\\
\bigwedge^{n-m-t}
G^*\ar[u]^{\theta}\ar[dr]_{\tilde\nu_\psi}\ar[rr]^{ x^*\wedge}  &&
\bigwedge^{n-t} G^*\ar[u]_{\theta}\ar[dl]^{\tilde\omega_{n-t}}\\
&(\bigwedge^t G^*)^*}.
$$ Using formula (\ref{antiderivation}) in section \ref{Preliminaries} one
sees easily that the first rectangle is commutative or
anticommutative. The second is commutative by
\ref{DUALIZATION_COMMUTATIVITY}. Finally
$$\tilde\nu_\psi(z^*)(y^*)=\tilde\omega( z^*\wedge  y^*\wedge x^*
)=\pm\tilde\omega(x^*\wedge z^*\wedge
y^*)=\pm\tilde\omega_{n-t}(x^*\wedge z^* )(y^*)$$ for all $z^*\in
\bigwedge^{n-m-t}G^*,\ y^*\in\bigwedge^{n-t} G^*$. An equally
simple computation shows that
$\tilde\omega_{n-t}\theta^{-1}\omega_t$ is the canonical map
$\bigwedge^t G\to (\bigwedge^tG^*)^*$.
\end{proof}

\begin{corollary}\label{IDENTIFICATION_2} If $G$ is free of rank
$n$, then there are (non-canonical) complex isomorphisms
$$
\mathcal C_\psi(t)\iso(\mathcal{C}_\psi(n-m-t))^*.
$$
\end{corollary}

\begin{proof} Since $\mathcal C^t(\psi)$ is dual to $\mathcal
C^{n-m-t}(\psi)$ (see [E], Chapter A.2.6.1, for example), the
corollary follows immediately from the preceding proposition.
\end{proof}
\bigskip

Now we shall associate a complex with the map $\phi:H\to G$,
generalizing the dual version of the classical Koszul complex.
Since the canonical map $\zeta_{H,G}:H^*\otimes G\to \Hom(H,G)$ is
an isomorphism, we may view $\phi$ as an element of the
$S(H^*)$-algebra
$$
S(H^*)\otimes \bigwedge G=\bigwedge_{S(H^*)}(S(H^*)\otimes G).
$$
We have $\phi^2=0$, so for every graded $S(H^*)$-module
$M=\oplus_{i}M_i$, the left multiplication by $\phi$ on
$M\otimes\bigwedge G $, denoted by $d_\phi$, gives rise to a
complex of $\bigwedge G$-modules
$$
\cdots\to M_{i-1}\otimes\bigwedge
G\stackrel{d_\phi}{\to}M_{i}\otimes\bigwedge
G\stackrel{d_\phi}{\to} M_{i+1}\otimes\bigwedge G\to\cdots.
$$
Two cases for $M$ lie in our interest: $M=S(H^*)$ with the natural
grading, and $M=D(H)$ where the structure of $D(H)$ as a graded
$S(H^*)$-module has been described in section \ref{Preliminaries}
and $M_i=D_{-i}(H)$ in this case. An easy calculation shows that
\begin{equation*}
d_\phi(h_1^{(k_1)}\ldots h_q^{(k_q)}\otimes y)=\sum_{j=1}^q
h^{(k_1)}_1 \ldots h^{(k_j-1)}_j \ldots h^{(k_q)}_q\otimes
\phi(h_j)\wedge y.
\end{equation*}
for all $h_1,\ldots,h_q\in H$, $y\in \bigwedge G$.

Let  $h_1,\ldots,h_l$ be a basis of $H$, $h_1^*,\ldots,h_l^*$ the
dual basis of $H^*$, and set $x=\phi(h_1)\wedge \ldots \wedge
\phi(h_l)\in \bigwedge G$. Then $\phi$, as an element of
$S(H^*)\otimes\bigwedge G$, has the presentation
\begin{equation*}
\phi=\sum_{j=1}^l h_j^*\otimes \phi(h_j), \end{equation*}
consequently $(1_{S(H^*)}\otimes x)\phi=\phi(1_{S(H^*)}\otimes
x)=0$ in $S(H^*)\otimes\bigwedge G$. So we may use the left
multiplication by $x$ on $\bigwedge G$, denoted by $\nu^{\phi}$,
in order to splice the complexes $(D(H)\otimes\bigwedge G,d_\phi)$
and $(S(H^*)\otimes\bigwedge G ,d_\phi)$ to a complex
$$
\cdots\to D_{1}(H)\otimes\bigwedge G
\stackrel{d_\phi}{\to}D_{0}(H)\otimes\bigwedge G
\stackrel{\nu^{\phi}}{\to} S_{0}(H^*)\otimes\bigwedge G
\stackrel{d_\phi}{\to}S_{1}(H^*)\otimes\bigwedge G \to\cdots.
$$
This complex, as a complex of $R$-modules, splits into direct
summands
\begin{align*}
\mathcal D_\phi(t):\quad 0\to D_{t}(H)\otimes \bigwedge^{0}
G\stackrel{d_\phi}{\to}\cdots \stackrel{d_\phi}{\to}
D_{0}(H)\otimes \bigwedge^{t} G & \stackrel{\nu^{\phi}}{\to}
S_{0}(H^*)\otimes \bigwedge^{t+l} G\stackrel{d_\phi}{\to}\cdots \\
& \stackrel{d_\phi}{\to}S_{p}(H^*)\otimes \bigwedge^{t+l+p} G
\to\cdots.
\end{align*}

Similarly as in the case of $\mathcal C_\psi(t)$ (see remark
\ref{USUALKOSZULINDENTIFICATION}), the connection map $\nu^{\phi}$
in $\mathcal D_\phi(t)$ depends on the basis taken for $H$.
Moreover there is an analogue with Proposition
\ref{IDENTIFICATION_1} and its proof.

\begin{proposition}\label{IDENTIFICATION_3}Suppose $G$ to be a free $R$-module of rank $n$.
Then  $\mathcal D_\phi(t)$ is (non-canonically) isomorphic to the
complex $\mathcal{D}^t(\phi)$ as defined in {\rm [BV1]}.
\end{proposition}

We are now able  to formulate a general result concerning the
connection of the generalized Koszul complex and its dual through
natural complex isomorphisms.
\begin{proposition}\label{DUALITY_EXTENDED}Let $\psi:G\to F$ be a homomorphism of
$R$-modules where $F$ is free of finite rank. Then the canonical
map $\theta:\bigwedge G^*\to (\bigwedge G)^*$ induces a natural
complex morphism
$$
\tau:\ \mathcal D_{\psi*}(t)\to\Big(\mathcal C_\psi(t)\Big)^*
$$
where the connection homomorphisms $\nu^{\psi^*}$ and $\nu_\psi$
are defined with respect to the same basis of $F$. If $\theta$ is
an isomorphism (for example, if $G$ is finitely generated and
projective),
 then $\tau$ is a complex isomorphism.
\end{proposition}

\begin{proof} We identify $D(F^*)= S(F)^*$ and $S(F)= S(F)^{**}$
by the corresponding graded algebra isomorphisms. Then we define
$\tau: D(F^*)\otimes \bigwedge G^*\to (\bigwedge G)^*\otimes
S(F)^*$ as the composition of the canonical map $D(F^*)\otimes
\bigwedge G^*\to \bigwedge G^*\otimes S(F)^*$ and $\theta\otimes
1_{S(F)^*}$. The map $\tau: S(F)\otimes \bigwedge G^*\to
(\bigwedge G)^*\otimes S(F)^{**}$ is defined analogously.

To show that $\tau$ is a complex morphism, we choose a basis
$f_1,\ldots,f_m$ for $F$. Let $f_1^*,\ldots,f_m^*$ be the dual
basis of $F^*$. Then $\psi^*=\sum_jf_j\otimes \psi^*(f_j^*)$ as an
element of $F\otimes G^*$, and
$$d_{\psi^*}=\sum_j f_j\, 1_A\otimes(\psi^*(f_j^*)\wedge\,)\quad\text{and}\quad
(\partial_\psi)^*=\sum_j\big(\leftharpoonup
\psi^*(f_j^*)\big)^*\otimes f_j\, 1_A$$ where $A$ stands for
$D(F^*),\ S(F),\ S(F)^*,\ S(F)^{**}$ respectively. Similarly
$$\nu^{\psi^*}=\psi^*(f_1^*)\wedge\ldots\wedge\psi^*(f_m^*)\wedge\quad\text{and}\quad
(\nu_\psi)^*=\big(\leftharpoonup\psi^*(f_1^*)\wedge\ldots\wedge\psi^*(f_m^*)\big)^*.$$
An application of Proposition \ref{DUALIZATION_COMMUTATIVITY} now
finishes the proof.

\end{proof}

\begin{remark} Suppose $G$ to be finitely generated free. Using the above result we get natural isomorphisms
$$
\mathcal C_\psi(t)\iso\Big(\mathcal
C_\psi(t)\Big)^{**}\iso\Big(\mathcal D_{\psi*}(t)\Big)^*.
$$
So we may canonically identify $\mathcal C_\psi(t)$ with the
generic resolution introduced by
 Buchsbaum and Eisenbud in [BE3].
\end{remark}

\section{Koszul Bicomplexes}\label{Koszul Bicomplexes}

Let $F$, $G$, $H$ be $R$-modules. We assume $H$ and $F$ to be free
of ranks $l$ and $m$. Let
$$
\CD H@>\phi>>G@>\psi>>F
\endCD
$$
be a complex. Due to Proposition \ref{DIAGRAMCOMMUTATIVITY} we can
assemble the  complexes $\mathcal C_\psi(t)$ and $\mathcal
D_\phi(t)$ to our main tool, the bicomplexes
$\mathcal{K}_{.,.}(t)$ \scriptsize
$$
\CD
\vdots &&\vdots &&\vdots &&\vdots\\
@VVV @VVV @VVV @VVV\\
\cdots H\otimes\bigwedge^{t+m}G\otimes
F^*@>>>\bigwedge^{t+m+1}G\otimes F^*@>\pm\nu^{\phi}>>
\bigwedge^{t+l+m+1}G\otimes F^*@>>>H^*\otimes\bigwedge^{t+l+m+2}G\otimes F^*\cdots\\
@VVV @V\partial_{\psi}VV @V\partial_{\psi}VV @VVV\\
\cdots
H\otimes\bigwedge^{t+m-1}G@>d_{\phi}>>\bigwedge^{t+m}G@>\pm\nu^{\phi}>>
\bigwedge^{t+l+m}G@>d_{\phi}>>H^*\otimes\bigwedge^{t+l+m+1}G\cdots\\
@V\pm\nu_\psi VV @V\pm\nu_\psi VV @V\pm\nu_\psi VV @V\pm\nu_\psi VV\\
\cdots H\otimes\bigwedge^{t-1}G @>d_{\phi}>>\bigwedge^{t}G
@>\pm\nu^{\phi}>>
\bigwedge^{t+l}G@>d_{\phi}>>H^*\otimes\bigwedge^{t+l+1}G\cdots\\
@VVV @V\partial_{\psi}VV @V\partial_{\psi}VV @VVV\\
\cdots H\otimes\bigwedge^{t-2}G\otimes
F@>>>\bigwedge^{t-1}G\otimes F@>\pm\nu^{\phi}>>
\bigwedge^{t+l-1}G\otimes F@>>>H^*\otimes\bigwedge^{t+l}G\otimes F\cdots\\
@VVV @VVV @VVV @VVV\\
\vdots &&\vdots &&\vdots &&\vdots\\
\endCD
$$
\normalsize The rows in the upper half are built from $\mathcal
D_\phi(t+m+j)$ tensored with $S_j(F)^*$, $j=0,1,\ldots$, while the
rows below are built from $\mathcal D_\phi(t-j)$ tensored with
$S_j(F)$, $j=0,1,\ldots$; we abbreviate $d_\phi\otimes 1_{S(F^*)}$
and $d_\phi\otimes 1_{S(F)}$ to $d_\phi$, and correspondingly
$\nu^\phi\otimes 1_{S(F^*)}$ and $\nu^\phi\otimes 1_{S(F)}$ to
$\nu^\phi$. The columns are obtained analogously: in western
direction we have to tensorize $D_i(H)$ with $\mathcal
C_\psi(t-i)$, $i=0,1,\ldots$, while going east we must tensorize
$S_i(H^*)$ with $\mathcal C_\psi(t+l+i)$, $i=0,1,\ldots$; as
before we shorten the complex maps to $\partial_\psi$ and
$\nu_\psi$. The signs of $\nu^\phi$ and $\nu_\psi$ are determined
by the formula in Proposition \ref{DIAGRAMCOMMUTATIVITY}.

To prove that the rectangles in $\mathcal{K}_{.,.}(t)$ are
anticommutative, is an easy exercise. We show for example that
$\partial_{\psi}d_{\phi}+d_{\phi}\partial_{\psi}=0$.

Let $M$ be a graded $S(H^*)$-module and $N$ a graded
$S(F)$-module. Then the diagram
$$
\CD
M\otimes\bigwedge G\otimes N@>d_{\phi}\otimes 1_N>>M\otimes\bigwedge G\otimes N\\
@V1_M\otimes \partial_{\psi}VV @VV1_M\otimes\partial_{\psi} V\\
M\otimes\bigwedge G\otimes N@>d_{\phi}\otimes 1_N>>M\otimes\bigwedge G\otimes N.\\
\endCD
$$
is anticommutative: as above let $h_1,\ldots,h_l$ and
$f_1,\ldots,f_m$ be bases of $H$ and $F$, and denote by
$h_1^*,\ldots,h_l^*$ and $f_1^*,\ldots,f_m^*$ the corresponding
dual bases. Then $\phi=\sum_ih_i^*\otimes\phi(h_i)$ and
$\psi=\sum_j\psi^*(f_j^*)\otimes f_j$ as elements of $H^*\otimes
G$ and $G^*\otimes F$ respectively. Then for $m\in M$, $y\in
\bigwedge G$, and $n\in N$
$$
(1_M\otimes\partial_{\psi})(d_{\phi}\otimes 1_N)(m\otimes y\otimes
n)=\sum_i\sum_j h^*_im\otimes\big((\phi(h_i)\wedge
y)\leftharpoonup\psi^*(f_j^*)\big)\otimes f_jn,
$$
while
$$
(d_{\phi}\otimes 1_N)(1_M\otimes\partial_{\psi})(m\otimes y\otimes
n)=\sum_j\sum_i h^*_im\otimes\big(\phi(h_i)\wedge
(y\leftharpoonup\psi^*(f_j^*))\big)\otimes f_jn.
$$
Since $\psi^*(f_j^*)(\phi(h_i))=f_j^*\psi\phi(h_i)=0$, the
associativity formula of Proposition \ref{DIAGRAMCOMMUTATIVITY}
yields the desired result.
\bigskip

\section{Grade Sensitivity}\label{Grade Sensitivity}

This section links up the study of the Koszul bicomplexes
introduced in the previous section with the study of the homology
of certain Koszul complexes. To avoid permanent repetitions we
shall assume throughout that the ring $R$ is {\it noetherian}.

For a linear map $\phi:H \to G$ of finitely generated free
$R$-modules we denote by $I_\phi$ the ideal in $R$ generated by
the maximal minors of a matrix representing $\phi$. We shall use
Macaulay's Theorem several times which says that $\grade I_\phi\le
|\rank G-\rank H|+1$ if $I_\phi\ne R$.

Our general assumption throughout the rest of this section will be
that {\it $H$, $G$ and $F$ are free $R$-modules of finite ranks
$l$, $n$ and $m$, and that
$$
\CD H@>\phi>>G@>\psi>>F
\endCD
$$
is a complex.}  Although much of what we will do, holds formally
for any $l$, $n$ and $m$, the applications will refer to the case
in which $n\ge m$ and $n\ge l$. So {\it we suppose that $r=n-m\ge
0$, $s=n-l\ge 0$.} We set $g=\grade I_\psi$ and $h=\grade I_\phi$.

A first question is which restrictions $g$ and $h$ are subjected
in a situation as the one pictured above. In the sequel we shall
give some answers to this question. The following result is a
simple consequence of Macaulay's Theorem.
\begin{proposition}\label{RESTRICTION INF} Set $\rho=r-l$.
\begin{enumerate}
\item[{\rm (1)}] If $g,h\ge 1$, then $\rho\ge 0$. \item[{\rm (2)}]
If $g>|\, \rho\, | +1$, then $I_\phi\subset \Rad I_\psi$, and in
particular $h\le g$. \item[{\rm (3)}] Moreover, if $g\ge r+1$,
then $I_\phi\subset I_\psi$.
\end{enumerate}
\end{proposition}

\begin{proof}If $g\ge 1$, then $M=\Coker \psi^*$ has $\rank r$. So $M^*=\Ker\psi$ has $\rank
r$, too. In the same way $h\ge 1$ implies that $\phi$ is
injective, so $\Im\phi$ has $\rank l$. Since
$\Ker\psi\supset\Im\phi$, we obtain the first part.

While proving (2) and (3) we may assume that $I_\psi\ne R$.
Suppose that $g>|\, \rho\, | +1$. Take a prime ideal $P\supset
I_\psi$ in $R$. Then $\grade(I_\psi R_P)\ge g$. If
$I_\phi\not\subset P$, then $(\Im\phi)_P$ would be a free direct
summand of $G_P$ of rank $l$, and therefore $\psi_P$ can be viewed
as a map from a free module of rank $ n-l$ to a free module of
rank $m$. So $\grade (I_\psi R_P)\le |\, \rho\, | +1$, which
contradicts the hypothesis. It follows that $I_\phi\subset P$
which implies that $I_\phi\subset \Rad I_\psi$.

In case $g=\grade I_\psi=r+1$, we consider a rank $1$ direct
summand $\tilde H$ of $H$. Let $\tilde\varphi$ be the restriction
of $\phi$ to $\tilde H$. From [BV4], Proposition 1, we draw that
$I_{\tilde\phi}\subset I_\psi$ in this case. Since $ I_\phi\subset
I_{\tilde\phi}$, the conclusion follows.
\end{proof}

Since $G$ is finitely generated, the generalized Koszul complexes
$\mathcal C_{\psi}(t)$ and $\mathcal D_{\phi}(t)$ (see section
\ref{Generalized Koszul Complexes}) have only a finite number of
non-vanishing components. To identify the homology, we fix their
graduations as follows: position 0 is held by the leftmost
non-zero module. The homology of $\mathcal C_{\psi}(t)$ and
$\mathcal D_{\phi}(t)$ behaves similarly as the homology of the
classical Koszul complex. The main result is the following.
\begin{theorem}\label{E-N HOMOLOGY}Set
$C=\Coker \psi$ and $D=\Coker \phi^*$. Furthermore let
$S_0(D)=R/I_\phi$, $S_{-1}(D)=\bigwedge^{s+1}\Coker\phi$,
$S_0(C)=R/I_\psi$ and $S_{-1}(C)=\bigwedge^{r+1}\Coker\psi^*$.
Then the following hold.
\begin{enumerate}
\item[\rm{(a)}]$H^i(\mathcal D_{\phi}(t))=0$ for $i<h$. Moreover,
if $t\le s+1$ and $\grade I_k(\phi)\ge n-k+1$ for all $k$ with
$l\ge k\ge 1$, then $\mathcal D_{\phi}(t)$ is a free resolution of
$S_{s-t}(D)$. (If $-1\le t\le s+1$, then it suffices to require
that $\grade I_\phi\ge s+1$.) \item[\rm{(b)}]$H^i(\mathcal
C_{\psi}(t))=0$ for $i<g$. Moreover, if $t\ge -1$ and $\grade
I_k(\psi)\ge n-k+1$ for all $k$ with $m\ge k\ge 1$, then $\mathcal
C_{\psi}(t)$ is a free resolution of $S_{t}(C)$. (If $-1\le t\le
r+1$, then it suffices to require that $\grade I_\psi\ge r+1$.)
\end{enumerate}
Finally, if $I_{\phi}=R$ ($I_{\psi}=R$), then all sequences
$\mathcal D_{\phi}(t)$ ($\mathcal C_{\psi}(t)$) are split exact.
\end{theorem}

\begin{proof}Proposition \ref{IDENTIFICATION_3}, Proposition \ref{DUALITY_EXTENDED} and
 Corollary  \ref{IDENTIFICATION_2}
 provide complex isomorphisms
$$
\mathcal{D}^t(\phi)\iso \mathcal D_{\phi}(t)\iso \mathcal
C_{\phi^*}(s-t),
$$
and similarly we obtain
$$
\mathcal{C}^t(\psi)\iso \mathcal C_{\psi}(t)\iso \mathcal
D_{\psi^*}(r-t).
$$
Since $\grade I_\phi=\grade I_{\phi^*}$ and $\grade I_\psi=\grade
I_{\psi^*}$, (a) follows from Proposition 2.1 in [BV1] while (b)
is obtained from Theorem A2.10,(c) in [E].
\end{proof}

By $\mathcal{C}_{.,.}(t)$ we shall denote the bicomplex which is
the lower part of the Koszul bicomplex in the previous section
(the rows below the second row). We rewrite this complex as

$$ \CD  0@>>> C_t^{0,0}
@>>>C_t^{1,0} @>d_\phi>>\cdots C_t^{t,0}
@>\pm\nu^{\phi}>>C_t^{t+1,0} @>d_\phi>>C_t^{t+2,0} && \cdots\\
&&@VVV @V\partial_\psi VV @VVV @V\partial_\psi VV @VVV \\
&&0@>>>C_t^{1,1} @>>>\cdots C_t^{t,1} @>\pm\nu^{\phi}>>
C_t^{t+1,1}@>>> C_t^{t+2,1} && \cdots\\ &&&&&& \vdots && \vdots && \vdots\\
&&&&&& C_t^{t,t}@>>> C_t^{t+1,t}@>>> C_t^{t+2,t}&& \cdots\\ &&&&&& @VVV @VVV @VVV\\
&&&&&& 0@>>> C_t^{t+1,t+1}@>>> C_t^{t+2,t+1} &&\cdots \\ &&&&&&
@VVV @VVV @VVV\\ &&&&&& \vdots && \vdots && \vdots
\endCD
$$
In other words,
$$C_t^{0,0}\ =\ \begin{cases}
D_t(H)\otimes\bigwedge^{0}G\otimes S_0(F)\quad &\text{if}\quad
0\le t,\\
S_0(H^*)\otimes\bigwedge^{t+l}G\otimes S_0(F)\quad
&\text{if}\quad -l\le t< 0,\\
S_{-t-l}(H^*)\otimes\bigwedge^{0}G\otimes S_0(F)\quad
&\text{if}\quad t< -l.\end{cases}
$$
The row homology of $\mathcal C_{.,.}(t)$ at $C_t^{p,q}$ is
denoted by $H_\phi^{p,q}$, the column homology by $H_\psi^{p,q}$.
Thus $H_\phi^{p,0}$ is the $p$-th homology module of $\mathcal
D_{\phi}(t)$.

Set $N^p=\Ker {(C_t^{p,0}\overset {\partial_\psi}\to C_t^{p,1})}$.
The canonical injections $N^p\to C_t^{p,0}$ yield a complex
morphism
$$ \CD 0@>>> N^0@>>>N^1 @>>>\cdots && N^p @>\bar{d_\phi}>>
N^{p+1}&&\cdots\\ &&
\parallel && @VVV && @VVV @VVV \\ 0@>>> C_t^{0,0}@>>>
C_t^{1,0}@>>>\cdots &&C_t^{p,0}@>d_\phi >> C_t^{p+1,0}&&\cdots
\endCD
$$ where the maps $\bar {d_\phi}$ are induced by $d_\phi$. The homology of
the first row $\mathcal N(t)$ at $N^p$ is denoted by $\bar{H}^p$.
We shall now investigate this homology.

For this purpose we extend $\mathcal{C}_{.,.}(t)$ to the complex
$\tilde{\mathcal{C}}_{.,.}(t)$ by setting $C_t^{p,-1}=N^p$. We
record some facts about the homology of
$\tilde{\mathcal{C}}_{.,.}(t)$. To avoid new symbols, the column
homology at $C_t^{p,q}$ is again denoted by $H_\psi^{p,q}$
(actually it differs from that of $\mathcal{C}_{.,.}(t)$ only at
$C^{p,0}$). We draw from Theorem \ref{E-N HOMOLOGY},(a) that
\begin{equation}
H_\phi^{p,q}=0\quad\text{ for $p<h$ and $q\ne-1$}.
\end{equation}
In case $r+1-g<p\le t$, we get from Theorem \ref{E-N HOMOLOGY},(b)
\begin{equation}
H_\psi^{p,q}=\begin{cases} \phantom{D_{t-p}(H)}0&\quad \text{for}\
0\le q<\min(p-(r+1-g),g),\\ D_{t-p}(H)\otimes S_p(C)&\quad
\text{for}\ q=p,\end{cases} \end{equation} and if $0\le t<p$, we
obtain
\begin{equation}
H_\psi^{p,q}=\begin{cases}
\phantom{S_{p-t-1}(H^*)}0&\quad\text{for}\ 0\le
q<\min(p-(\rho+2-g),g),\\S_{p-t-1}(H^*)\otimes S_{p+l-1}(C)&\quad
\text{for}\ q=p+l-1.\end{cases}
\end{equation}
For $q\ge -1$ we consider the $q$th row $\tilde C_{.,q}$ of
$\tilde C_{.,.}$ and its image complex $\partial_\psi (\tilde
C_{.,q})$ in $\tilde C_{.,q+1}$. We set $E^{p,q}=H^p(\partial_\psi
(\tilde C_{.,q-1}))$ for $q\ge 0$. Then there are exact sequences
\begin{equation}
H_\phi^{i-(j+1),j} \to E^{i-(j+1),j+1}\to E^{i-j,j} \to
H_\phi^{i-j,j}
\end{equation}
if $H_\psi^{i-(j+1),j}=H_\psi^{i-j,j}=0$, and because of (4) and
(5) this holds if
\begin{align*}
r+1-g<i-(j+1)\le t\ \text{and}\
&j<\min\big(i-(j+1)-(r+1-g),g\big)\ \text{or}\\
0\le t<i-(j+1)\ \text{and}\ &j<\min\big(i-(j+1)-(\rho+2-g),g\big).
\end{align*}

\begin{theorem}\label{FUNDAMENTAL} Let $t\ge 0$.
Set $C=\Coker \psi$, $\rho=r-l$. Then $\bar{H}^i= 0$ for
$i=0,\ldots,\min(2,h-1)$. Suppose that $g>|\, \rho\, |+1$.
\begin{enumerate}
\item[\rm{(i)}] If $3\le 2t+1\le h$, then there is an injection $
D_{0}(H) \otimes S_{t} (C)\to \bar{H}^{2t+1}$ which is an
isomorphism if $2t+1<h$; \item[\rm{(ii)}] $\bar{H}^{i}=0$ for
$2t+2\le i< \min(h,2t+g-|\, \rho\, |+1)$; \item[\rm{(iii)}] if
 $2t+g-\mid\rho\mid+1<h$, then
$\bar{H}^{2t+g-\mid\rho\mid+1}=H_\psi^{t+1,\ t+g-\mid\rho\mid-1}$.
\end{enumerate}
\end{theorem}

\begin{proof}
We may obviously assume that $h>0$. From Proposition
\ref{RESTRICTION INF} we then obtain that $\rho\ge 0$. Set
$\mu=\min(2,h-1)$. Consider the diagram
$$
\CD  0@>>> N^0@>>>N^1 @>>>N^2 @>\bar d_{\phi}>> N^{3}\\ &&
\parallel && @VVV @VVV @VVV
\\ 0@>>> C_t^{0,0}@>>> C_t^{1,0}@>>>C_t^{2,0}@>d_\phi >> C_t^{3,0}\\ &&
@VVV  @VVV  @VVV @VVV \\
0@>>>0@>>> \Im \partial_\psi^{1,0}@>>> \Im
\partial_\psi^{2,0} @>>> \Im \partial_\psi^{3,0}\\
 \endCD
$$
The middle row has trivial homology at $C_t^{p,0}$ for $p\le\mu$
(Theorem \ref{E-N HOMOLOGY}). Because of $h\ge 1$, the row
homology at $\Im
\partial_\psi^{1,0}$ is zero since the homomorphism $\Im
\partial_\psi^{1,0}\to \Im \partial_\psi^{2,0}$ is the restriction of the
injective homomorphism $ C_t^{1,1}\overset{d_\phi}\to C_t^{2,1}$.
Now we use the long exact homology sequence to get the first
statement.

Next let $g>\rho+1$ and $3\le 2t+1\le h$. Remark that $h\le g$ by
Proposition \ref{RESTRICTION INF},(2). Since
$H_\psi^{2t-j,j}=H_\psi^{2t+1-j,j}=0$ for $j=0,\ldots,t-1$, we
obtain exact sequences
$$ H_\phi^{2t-j,j} \to E^{2t-j,j+1}\to E^{2t+1-j,j} \to H_\phi^{2t+1-j,j}
$$ for these $j$. Furthermore
$2t+1\le h$ by assumption. Therefore we obtain the ``southwest''
isomorphisms
 $$E^{2t,1}\iso\dots\iso
E^{t+2,t-1}\iso E^{t+1,t} $$ and an injection $E^{2t,1}\to
E^{2t+1,0}=\bar H^{2t+1}$ which is an isomorphism if even $2t+1<
h$.
 In the diagram\scriptsize
$$\CD  &&
\Im\partial_\psi^{t,t-1}@>>>
\Im\partial_\psi^{t+1,t-1}=\Ker\partial_\psi^{t+1,t}@>>>
\Im\partial_\psi^{t+2,t-1}=\Ker\partial_\psi^{t+2,t}\\
&& @VVV @VVV @VVV  \\
0@>>> C_t^{t,t}@>>>
C_t^{t+1,t}@>>> C_t^{t+2,t}\\
&&@VVV @VVV @VVV  \\ 0 @>>>D_{0}(H)\otimes S_{t}(C)@>>> \Im
\partial_\psi^{t+1,t}@>d_\phi>>
\Im\partial_\psi^{t+2,t}
\endCD
$$\normalsize
we abbreviate
$\partial_\psi^{p,q}=(C_t^{p,q}\overset{\partial_\psi}\to
C_t^{p,q+1}$) and $d_\phi^{p,q}=(C_t^{p,q}\overset{d_\phi}\to
C_t^{p+1,q})$. It is induced by $\mathcal{C}_{.,.}(t)$ and has
exact middle row. Since $d_\phi^{t+1,t+1}$ is injective, we obtain
$D_{0}(H) \otimes S_{t} (C)\iso E^{t+1,t}$. This together with the
considerations in the previous paragraph proves (i).

To prove (ii), let $2t+2\le i< \min(h,2t+g-\rho+1)$. Going
``southwest'' once more, we obtain isomorphisms
$$ \bar H^i=E^{i,0}\iso
E^{i-1,1}\iso\dots\iso E^{t+1,i-t-1}.$$ Since $d_\phi^{t+1,i-t-1}$
is injective, the claim follows.

In case $2t+g-\rho+1<h$, we get
$$ \bar H^{2t+g-\rho+1}=E^{2t+g-\rho+1,0}\iso
E^{2t+g-\rho,1}\iso\dots\iso E^{t+2,t+g-\rho-1}.$$ Furthermore the
diagram
$$ \CD
&&\Im\partial_\psi^{t+1,t+g-\rho-2}
@>>>\Im\partial_\psi^{t+2,t+g-\rho-2}
@>>>\Im\partial_\psi^{t+3,t+g-\rho-2}\\
&& @VVV @VVV @VVV \\
0@>>>\Ker\partial_\psi^{t+1,t+g-\rho-1}
@>d_\phi>>\Ker\partial_\psi^{t+2,t+g-\rho-1}
@>>>\Ker\partial_\psi^{t+3,t+g-\rho-1}\\
&& @VVV @VVV @VVV
\\  0@>>> H_\psi^{t+1,t+g-\rho-1}@>>>0@>>>0
\endCD
$$
has row homology zero at $\Ker\partial_\psi^{t+1,t+g-\rho-1}$
since $d_\phi^{t+1,t+g-\rho-1}$ is injective. If we can show that
the row homology at $\Ker\partial_\psi^{t+2,t+g-\rho-1}$ also
vanishes, we shall obtain (iii). The diagram
$$
\CD 0@>>>\Ker\partial_\psi^{t+1,t+g-\rho-1}
@>>>\Ker\partial_\psi^{t+2,t+g-\rho-1}
@>>>\Ker\partial_\psi^{t+3,t+g-\rho-1}\\
&& @VVV @VVV @VVV \\
0@>>> C^{t+1,t+g-\rho-1}@>>> C^{t+2,t+g-\rho-1}@>>>
C^{t+3,t+g-\rho-1}\\
&& @VVV @VVV @VVV\\
0@>>>\Im\partial_\psi^{t+1,t+g-\rho-1}
@>d_\phi>>\Im\partial_\psi^{t+2,t+g-\rho-1}
@>>>\Im\partial_\psi^{t+3,t+g-\rho-1}\\
\endCD
$$
has exact middle row. Since $d_\phi^{t+1,t+g-\rho}$ is injective,
we get the desired result.
\end{proof}

Provided $g\ge r+1$, we can compute the homology $\bar H^i$ of
$\mathcal N^{^.}(t)$ outside the gap $2t+1\le i<\min (h,2t+3+l)$.
\begin{theorem}\label{FUNDAMENTAL_2} As
before let $t\ge 0$, and
 set $C=\Coker \psi$, $\rho=r-l$. Suppose that $g\ge r+1$.
 Then
\begin{enumerate}
\item[\rm{(i)}] $$\bar{H}^{i}=\begin{cases}
D_{t-\frac{i-1}2}(H)\otimes S_{\frac{i-1}2}(C)\quad &\text{if}\
3\le i<\min(h,2t+3),\ i\not\equiv 0\ (2),\\
\phantom{D_{t-\frac{i-1}2}(H)}0\quad &\text{if}\ 3\le
i<\min(h,2t+3),\ i\equiv 0\ (2);\end{cases}$$

\item[\rm{(ii)}]$$\bar{H}^{i}=\begin{cases}S_{\frac{i-l}2-t-1}(H)\otimes
S_{\frac{i+l}2-1}(C)\quad &\text{if}\ \ 2(t+1)+l\le i<h,\
i\equiv l\ (2),\\
\phantom{D_{t-\frac{i-1}2}(H)}0\quad &\text{if}\ \ 2(t+1)+l\le
i<h,\ i\not\equiv l\ (2).\end{cases}$$
\end{enumerate}
\end{theorem}
\begin{proof} For (i) let $3\le i< \min(h,2t+3)$. Then
\begin{equation*}
H_\psi^{i-(j+1),j}=H_\psi^{i-j,j}=0\quad\text{for}\quad\begin{cases}\
j=0,\ldots,\frac{i-2}2&\ \text{if}\ i\equiv 0\ (2)\\
\ j=0,\ldots,\frac{i-3}2&\ \text{if}\ i\not\equiv 0\ (2)
\end{cases}
\end{equation*}
Since $H_\phi^{i-(j+1),j}=H_\phi^{i-j,j}=0$ because of $i<h$, the
sequences (6) may be used to obtain the isomorphisms
$$ \bar
H^{i}=E^{i,0}\iso E^{i-1,1}\iso\dots\iso E^{\frac {i}2,\frac
{i}2}.
$$ if $i$ is even, and $$ \bar H^i=E^{i,0}\iso
E^{i-1,1}\iso\dots\iso E^{\frac {i+1}2,\frac {i-1}2}
$$ if $i$ is odd.

In the last case we substitute $t$ by $\frac{i-1}2$ and
$D_0(H)\otimes S_t(C)$ by $D_{t-\frac{i-1}2}\otimes
S_{\frac{i-1}2}(C)$ in the second diagram in the proof of the
previous theorem. The resulting diagram has the same properties as
the original one. Therefore $$D_{t-\frac{i-1}2}\otimes
S_{\frac{i-1}2}(C)\iso E^{\frac {i+1}2,\frac {i-1}2}.$$ If $i$ is
even, we look at the sequence
\begin{equation*}
\Im\partial_\psi^{\frac{i-2}2,\frac{i-2}2}\to
\Im\partial_\psi^{\frac i2,\frac{i-2}2}\to \Im\partial_\psi^{
\frac{i+2}2,\frac {i-2}2}.
\end{equation*}
But $\partial_\psi^{\frac{i-2}2,\frac{i-2}2}=0$, and the second
map is induced by the injective map $d_\phi:C^{\frac i2,\frac
i2}\to C^{\frac {i+2}2,\frac i2}$. Therefore $E^{\frac {i}2,\frac
{i}2}=0$.

To prove (ii) one uses the same strategy.
\end{proof}

\begin{remarks}
The proofs of Theorems \ref{FUNDAMENTAL} and \ref{FUNDAMENTAL_2}
partially consist of a repetition of arguments used in the proof
of Proposition 1 in [BV4].
\end{remarks}

One can extend the assertions of \ref{FUNDAMENTAL} and
\ref{FUNDAMENTAL_2} to negative integers $t$. We touch briefly
upon this case omitting the very similar proofs.
\begin{theorem}\label{FUNDAMENTALnegativ}Let $t<0$ be an integer.
As above we set $C=\Coker \psi$, $\rho=r-l$.
\begin{enumerate}
\item[\rm{(a)}]Suppose that $g>|\, \rho\, |+1$. Then
\begin{enumerate}
\item[\rm{(i)}]$\bar{H}^{i}=0$ for $0\le i< \min(h,\max(2,t+g-|\,
\rho\, |))$ and \item[\rm{(ii)}]$\bar{H}^{t+g-|\, \rho\,
|}=H_\psi^{0,t+g-|\, \rho\, |-1}$ if $2\le t+g-|\, \rho\, |<h$.
\end{enumerate} \item[\rm{(b)}]Suppose that $g\ge r+1$.
\begin{enumerate}\item[\rm{(i)}] If $-l<t$, then
\begin{equation*}
\bar H^{i}=\begin{cases} S_{\frac{i-l-t-1}2}(H^*)\otimes
S_{\frac{i+l+t-1}2}(C)&\text{if}\ l+t+1\le i<h,
i+l+t\not\equiv 0\ (2),\\
\phantom{S_{\frac{i-l-t-1}2}(H^*)}0&\text{otherwise if}\ 0\le
i<h.\end{cases}\end{equation*} \item[\rm{(ii)}] If $t\le -l$, then
$\bar H^{i}= 0$ for $i=0,\ldots,\min(2,h-1)$ and
\begin{equation*}\bar
H^{i}=\begin{cases}S_{\frac{i-1}2-t-l}(H^*)\otimes
S_{\frac{i-1}2}(C)&\quad \text{if}\ 3\le i<h,\ i\not\equiv 0\
(2),\\
\phantom{S_{\frac{i-1}2-t-l}(H^*)}&\quad\text{if}\ 3\le i<h,\
i\equiv 0\ (2).\end{cases}\end{equation*}
\end{enumerate}
\end{enumerate}
\end{theorem}

\begin{remark} \label{N identification1} In view of what we will do in the next section, we record here
another interpretation of the complex $\mathcal N(t)$.

With the notion introduced above set $M=\Coker\psi^*$. By
$\bar{\lambda}:M\to H^*$ we denote the linear map induced by
$\phi^*$. Then there is a canonical complex isomorphism
$$
\CD \mathcal N(t)@>>>\Big(\mathcal C_{\bar{\lambda}}(t)\Big)^*.
\endCD
$$
To see this we substitute $\phi^*$ for $\psi$ in Proposition
\ref{DUALITY_EXTENDED}. We then obtain a natural complex
isomorphism
$$
\CD \mathcal D_{\phi}(t)@>\tau>>\Big(\mathcal
C_{\phi^*}(t)\Big)^*.
\endCD
$$
Obviously $\Big(\mathcal C_{\bar{\lambda}}(t)\Big)^*$ may be
viewed as a subcomplex of $\Big(\mathcal C_{\phi^*}(t)\Big)^*$,
and
$$\CD\tau(\mathcal N(t))=\Big(\mathcal C_{\bar{\lambda}}(t)\Big)^*,\endCD$$
 since
  $$
\Ker\big((\bigwedge^p G^*)^*\overset{(d_{\psi^*})^*}\to(
F^*\otimes\bigwedge^{p-1} G^*)^*\big)\iso (\bigwedge^p M)^*.
 $$
\end{remark}

\section{Koszul complexes in projective dimension one} We shall
now investigate the homology of the generalized Koszul complexes
in free (projective) dimension $1$.

In this section {\it $R$ is a noetherian ring, and $M$ an
$R$-module which has a presentation
$$
\CD 0@>>>\mathcal{F}@>\chi>>\mathcal{G}@>>>M@>>>0
\endCD
$$
where $\mathcal{F}$, $\mathcal{G}$ are free modules of ranks $m$
and $n$.} Then in particular $r=n-m\ge 0$.

In the sequel we consider $R$-homomorphisms
$\bar\lambda:M\to\mathcal{H}$ into a finite free $R$-module
$\mathcal{H}$ of rank $l\le n$. By $\lambda$ we denote the
corresponding lifted maps $\mathcal{G}\to \mathcal{H}$. Dualizing
$\mathcal F\overset {\chi}\to\mathcal G \overset{\lambda}\to
\mathcal H$ we go back to the situation previously studied. So we
set $F=\mathcal{F}^*$, $G=\mathcal{G}^*$, $H=\mathcal{H}^*$,
$\psi=\chi^*$, $\phi=\lambda^*$, and $C=\Coker \psi$.
\smallskip

We say that a homomorphism of finite free $R$-modules is {\it
minimal} if the entries of a representing matrix generate a proper
ideal of $R$.

\begin{theorem}\label{EXTENSION H-B submaximal}Set $\rho=r-l$ and $k=r+1-g$.
Suppose that $g=\grade I_\chi>|\,\rho\,|+1$.
\begin{enumerate}
\item[\rm{(a)}] $I_\lambda\subset \Rad I_\chi$, and in particular
$\grade I_\lambda\le g$. \item[\rm{(b)}] The following conditions
are equivalent:
\begin{enumerate}
\item[\rm{(1)}] $\grade I_\lambda>|\,\rho\,|+1$; \item[\rm{(2)}]
$\Rad I_\chi=\Rad I_\lambda$.
\end{enumerate}
\item[\rm{(c)}] Suppose that $I_\chi\ne R$ and that there is a
$\bar\lambda$ such that $\grade I_\lambda>|\,\rho\,|+1$. Then
$l=k+1$, $r\ge l$, and $r-k$ is odd.
\begin{enumerate}\item
[$\rm{(c_1)}$] Let $r=k+1$ $(=l)$. The sequence
$0\to\mathcal{F}\overset{\chi}\to\mathcal{G}\overset{\lambda}\to\mathcal
H$ as well as its dual $0\to H\overset{\phi}\to G\overset{\psi}\to
F$ are exact. This can occur only if $I_\chi=I_\lambda$.
Furthermore $m=1$ occurs if and only if the $i$th entry of a
matrix for $\chi$ is $(-1)^i$ times the minor of a matrix for
$\lambda$ by cancelling the $i$th row. \item [$\rm{(c_2)}$] Let
$r\ge k+3$ (i. e. $r\ge l+2)$. If $\chi$ is minimal, then
 $m\le k+1$ $(= l)$. If $\lambda$ is minimal, then
 $m> k$ $(= l-1)$.
\end{enumerate}
\end{enumerate}
\end{theorem}

\begin{proof}(a) With respect to the assumption, Proposition \ref{RESTRICTION
INF},(2) yields $I_\lambda= I_{\lambda^*}\subset\Rad
I_{\chi^*}=\Rad I_{\chi}$.
\smallskip

(b), (1)$\Rightarrow$(2) is an immediate consequence of
Proposition \ref{RESTRICTION INF},(2) while (2)$\Rightarrow$(1) is
trivial.
\smallskip

Next we prove the main statement of (c). From (b) we draw that
$\Rad I_\chi=\Rad I_\lambda$. Since $h:=\grade I_\lambda=g>1$,
Proposition \ref{RESTRICTION INF},(1) implies that $\rho\ge 0$ so
$r\ge l$, and $r-k+1=g>\rho+1=r-l+1\ge 1$ implies that $l>k$ and
$r-k>0$.

Assume that $r-k$ is even. Then $r-k\ge 2$. Consider the complex
$\mathcal N(\frac {r-k}2)$ defined in the previous section, and
observe that
$$
N^{r-k+1}(\frac{r-k}2)=S_{\frac{r-k}2}(H^*)\otimes(\bigwedge^{r+l-k}
M)^*
$$
since
$$
N^{r-k+1}(\frac{r-k}2)=\Ker\big(S_{\frac{r-k}2}(H^*)\otimes\bigwedge^{r+l-k}
G\to S_{\frac{r-k}2}(H^*)\otimes\bigwedge^{r+l-k-1} G\otimes
F\big).
$$
Now $l>k$ implies that $(\bigwedge^{r+l-k} M)^*=0$. So
$N^{r-k+1}(\frac{r-k}2)=0$. It follows that
$\bar{H}^{r-k+1}(\frac{r-k}2)=0$. On the other hand, using Theorem
\ref{FUNDAMENTAL} (a) we obtain an exact sequence
$$
\CD 0@>>>D_{0}(H)\otimes
S_{\frac{r-k}2}(C)@>>>\bar{H}^{r-k+1}(\frac{r-k}2).
\endCD
$$
Because $\bar{H}^{r-k+1}(\frac{r-k}2)=0$, we get $D_{0}(H)\otimes
S_{\frac{r-k}2}(C)=0$ and consequently $C=0$ which is in
contradiction with $I_\chi\ne R$. So $r-k$ must be odd.

If $r-k=1$, then $l=k+1$ since $l\le r$. Suppose that $r-k\ge 3$
(and odd). Then
$$
N^{r-k}(\frac{r-k-1}2)=
S_{\frac{r-k-1}2}(H^*)\otimes(\bigwedge^{r+l-k-1} M)^*.
$$
On the other hand we draw from Theorem \ref{FUNDAMENTAL},(a) that
$$
\bar{H}^{r-k}(\frac {r-k-1}2)= D_{0}(H)\otimes
S_{\frac{r-k-1}2}(C).
$$
If $l>k+1$, then $(\bigwedge^{r+l-k-1} M)^*=0$, so
$S_{\frac{r-k-1}2}(C)=0$, a contradiction. It follows that
$l=k+1$.

$\rm{(c_1)}$ The first statement is an immediate consequence of
the Buchsbaum-Eisenbud acyclicity criterion (see [E], Theorem
20.9). The second statement is also due to Buchsbaum and Eisenbud
(see [N], Ch. 7, Theorem 3, or Corollary 5.1 in [BE1]). (Of course
the third statement is a special case of the Theorem of
Hilbert-Burch.)

To prove the first claim of $\rm{(c_2)}$ let $\chi$ be minimal.
Localize at a prime ideal $P$ which contains $I_1(\chi)$. Then
$g\le \grade I_\chi R_P$. But $g<\grade I_\chi R_P$ is impossible:
otherwise, since $I_\lambda\subset\Rad I_\chi\subset P$ in view of
(a), we have $\grade I_\lambda R_P> |\, \rho\, |+1$ and
consequently $l=r+1-\grade I_\chi R_P+1\le k$ in contradiction
with the first claim under (c). So we may assume $R$ to be local.
Since $r\ge k+3$, $l=k+1$ and $M$ is free in depth $1$, we get
$N^{r-k}(\frac{r-k-1}2)= S_{\frac{r-k-1}2}(H^*)$ and
$N^{r-k+1}(\frac{r-k-1}2)=0$. So $\bar
H^{r-k}(\frac{r-k-1}2)=S_{\frac {r-k-1}2}(C)$ is a quotient of
$S_{\frac{r-k-1}2}(H^*)$. Therefore the minimal number of
generators of $C$ cannot be greater than the minimal number of
generators of $H$. It follows that $m\le l=k+1$.

Since $r\ge k+3$ and $l=k+1$, we deduce that $\rho\ge 2$. So
$h=g\ge 4$. To prove the second statement of $\rm{(c_2)}$ we
dualize the sequence $0\to \mathcal F\overset\chi\to \mathcal
G\overset\lambda\to\mathcal H$. Set $r'=n-l$ and $k'=r'+1-h$. From
the first claim under (c) we draw that $m=k'+1$. Since $h\ge 4$,
we get $r'\ge k'+3$.
 Now let
$\lambda$ be minimal.  Applying the first part of $\rm(c_2)$, we
obtain $(k+1=)\ l\le k'+1=m$.
\end{proof}

\begin{corollary}\label{h_limit} Set $\rho=r-l$ and suppose that
$\infty>\grade I_\chi>|\rho |+1$. Then $\grade I_\lambda\le |\,
\rho\, |+2$ for every $\bar\lambda$. Moreover, if $\grade
I_\lambda=|\, \rho\, |+2$, then $\grade I_\chi$, $\grade
I_\lambda$, and $\rho$  are even.
\end{corollary}

\begin{proof}If $\infty>\grade
I_\lambda>|\, \rho\, |+1$, then the above theorem implies that
$\grade I_\lambda=\grade I_\chi= r+1-k=r-l+2=|\, \rho\, |+2$.
Since $r-k$ is odd, the integers $\grade I_\chi,\grade
I_\lambda,\rho$ must be even.
\end{proof}

\begin{remark}
Let $k\in \mathbb{N}$. If $l=k+1$ and $r-k>0$ is odd, then, in the
cases listed under (c), there are always maps $\chi:R^m\to R^n$
and $\lambda:R^n\to R^l$ such that $\lambda \chi=0$ and $\grade
I_\chi=\grade I_\lambda=r-k+1$, provided there is a regular
sequence of length $r-k+1$ in $R$. (See [I], section 3.2 for an
example.)
\end{remark}

The following criterion could occasionally be useful.
\begin{corollary}\label{MatrixCriterion}Let $l\le n,\ m\le n$ and set $\rho=n-m-l$. Furthermore
let $\mathcal{A}$ be an $l\times n$-matrix and $\mathcal{B}$ be an
$n\times m$-matrix with entries in $R$. Set $h=\grade
I_\mathcal{A}$,
 $g=\grade I_\mathcal{B}$, and suppose that $\infty>h,g>|\, \rho\, |+1$. Then
$\mathcal{A}\mathcal{B}\neq 0$ in each of the following cases:
\begin{enumerate}
\item[\rm{(1)}]$h\neq g$; \item[\rm{(2)}]$\rho$ odd;
\item[\rm{(3)}]$g\neq |\, \rho\, |+2$; \item[\rm{(4)}]the ideals
generated by the entries of $\mathcal{A}$ and $\mathcal{B}$ are
proper ideals of $R$ and $l\neq m$ and $l\neq n-m$.
\end{enumerate}
\end{corollary}

In the second part of this section we treat the case in which {\it
$\grade I_{\chi}$ has the greatest possible value $n-m+1$}.
\smallskip

The following result may be seen as an extension of the
Hilbert-Burch Theorem.
\begin{theorem}\label{EXTENSION H-B}Suppose that
$\grade I_\chi=r+1$. Then $I_\lambda\subset I_\chi$, and in
particular $\grade I_\lambda\le r+1$. Set $\rho=r-l$.
\begin{enumerate}
\item[{\rm(a)}] If there is a $\bar\lambda$ such that $\grade
I_\lambda>|\, \rho\, |+1$, then $l=1$ and $r$ is odd.
\item[{\rm(b)}]Suppose in addition that $\chi$ is minimal. Then
the following are equivalent:
\begin{enumerate}
\item[{\rm(1)}] There is a $\bar\lambda$ such that $\grade
I_\lambda>|\, \rho\, |+1$; \item[{\rm(2)}] $l=1$ and {\rm (i)}
$r=1$ or {\rm (ii)} $m=1$ and $r\ge 3$ is odd.
\end{enumerate}
\item[{\rm(c)}]The following conditions are equivalent:
\begin{enumerate}
\item[{\rm(1')}] $\grade I_\lambda>|\, \rho\, |+1$;
\item[{\rm(2')}] $I_\lambda=I_\chi$.
\end{enumerate}
\end{enumerate}
\end{theorem}

\begin{proof} By Proposition \ref{RESTRICTION INF},(3)
we obtain $ I_\lambda= I_{\lambda^*}\subset I_{\chi^*}=I_\chi$.
Since $I_\chi\ne R$ by assumption, the first part of the theorem
is clear.

Theorem \ref{EXTENSION H-B submaximal},(c) implies (a) and also
the implication $(1) \Rightarrow (2)$ of (b) since $k=0$. The
other direction of (b) is a simple exercise.

It remains to prove (c), (1') $\Rightarrow$ (2') (since (2')
$\Rightarrow $ (1') is trivial): We may assume $R$ to be local,
and, using the uniqueness of minimal free resolutions, we can
easily reduce to the case in which $\chi$ is minimal. According to
what we have proved already, it follows that $l=1$ and (i) $r=1$
or (ii) $m=1$ and $r\ge 3$ is odd. If $r=1$, then we can apply the
Hilbert-Burch Theorem to get the desired equality
$I_\lambda=I_\chi$. If $m=1$ we look at the exact sequence
$$
\CD R@>\lambda^*>>\mathcal{G}^*@>\chi^*>>R.
\endCD
$$
which satisfies the hypothesis of Proposition \ref{RESTRICTION
INF},(3) since $\grade I_\lambda >n-m$ by assumption. So
$$
I_\chi=  I_{\chi^*}\subset I_{\lambda^*}=I_\lambda.
$$
\end{proof}

\begin{corollary}\label{EXTENSION H-B 2}Suppose that $\grade I_\chi=r+1$. Then the following conditions are equivalent:
\begin{enumerate}
\item[{\rm(1)}] there is a $\bar\lambda$ with $\grade
I_\lambda=n-l+1$; \item[{\rm(2)}] $l=1$,  $m=1$ and $r\ge 1$ is
odd.
\end{enumerate}
\end{corollary}

\begin{proof} Only (1) $\Rightarrow$ (2) requires a proof.
First $h,g\ge 1$, so $\rho\ge 0$. Then we have $n-l+1>\rho+1$.
Consequently $I_\lambda=I_\chi$. From Theorem \ref{EXTENSION H-B},
(a) we obtain that $l=1$ and $r\ge 1$ is odd. Since $r=n-l$, it
follows that $m=1$.
\end{proof}

\begin{corollary}\label{EXTENSION H-B 3}Suppose $R$ to be local and $\grade I_\chi=r+1$.
Then the following conditions are equivalent:
\begin{enumerate}
\item[\rm{(1)}] there is a $\bar\lambda$ with $\grade
I_\lambda>|\, \rho\, |+1$; \item[\rm{(2)}] $l=1$ and {\rm (i)}
$r=1$ or {\rm (ii)} $M$ has a minimal resolution
$$
\CD 0@>>>R@>>>R^{2k}@>>>M@>>>0
\endCD
$$
where $k\ge 2$.
\end{enumerate}
\end{corollary}
\bigskip
Our next goal is the description of the homology
of
\begin{align*}
\mathcal C_{\bar{\lambda}}(t):\cdots\to \bigwedge^{t+l+p} M\otimes
S_{p}(\mathcal{H})^*
 \stackrel{\partial_{\bar{\lambda}}}{\to}\cdots\stackrel{\partial_{\bar{\lambda}}}{\to}
 \bigwedge^{t+l} M &\otimes S_{0}(\mathcal{H})^*
\stackrel{\nu_{\bar{\lambda}}}{\to}\bigwedge^{t} M\otimes
S_{0}(\mathcal{H})
 \stackrel{\partial_{\bar{\lambda}}}{\to} \\
& \cdots \stackrel{\partial_{\bar{\lambda}}}{\to} \bigwedge^{0} M
\otimes S_{t}(\mathcal{H})\to 0.
\end{align*}
For this purpose we refer to the the upper part of the Koszul
bicomplex $\mathcal K_{.,.}(t)$ of section \ref{Koszul
Bicomplexes} (the rows above the second row) which we denote by
$\mathcal{B}_{.,.}(t)$. We rewrite this complex as
$$
\CD
 &&\vdots &&&&\vdots && \vdots &&\\
 &&@VVV  &&@VVV @VVV  \\
\cdots@>>> B_t^{0,-1}@>>>
\cdots @>d_\phi>>B_t^{t,-1} @>\pm\nu^{\phi}>>B_t^{t+1,-1}@>d_\phi>>\cdots\\
&&@VVV && @V\partial_\psi VV  @V\partial_\psi VV   \\
\cdots@>>>B_t^{0,0}@>>> \cdots@>>> B_t^{t,0} @>\pm\nu^{\phi}>>
B_t^{t+1,0}@>>>\cdots
\endCD
$$
where
 $$
 B^{0,0}_t=\begin{cases}
D_t(H)\otimes\bigwedge^{m}G\otimes S_0(F)^*&\quad\text{ if $0\le
t$},\\
S_0(H^*)\otimes\bigwedge^{t+l+m}G\otimes S_0(F)^*&\quad\text{ if
$-l\le t< 0$},\\
S_{-t-l}(H^*)\otimes\bigwedge^{m}G\otimes S_0(F)^*&\quad\text{ if
$t< -l$.}
\end{cases}$$
Set $M^p=\Coker
(B_t^{p,-1}\overset{\partial_\psi}\to B_t^{p,0})$. The canonical
surjection $B_t^{p,0}\to M^p$ yields a complex morphism
$$
\CD \cdots @>>> B_t^{-1,0}@>>> B_t^{0,0}@>>>\cdots @>d_\phi
>>B_t^{p,0}@>>>
B_t^{p+1,0}@>>>\cdots\\
 && @VVV @VVV &&
@VVV  @VVV\\
\cdots @>>>M^{-1}@>>>M^0 @>>>\cdots
@>\bar{d_\phi}>>M^p @>>> M^{p+1}@>>>\cdots\\
\endCD
$$ where the maps $\bar {d_\phi}$ are induced by $d_\phi$. The lower row is denoted by
$\mathcal M(t)$.

We obtain an analogue with Remark \ref{N identification1}.
\begin{proposition}\label{N identification2}Set $\rho=n-m-l$ as before.
Then there is a (non-canonical) complex isomorphism
$$
\CD \mathcal M(t)@>>>\mathcal C_{\bar{\lambda}}(\rho-t)
\endCD
$$
\end{proposition}

\begin{proof}Corollary  \ref{IDENTIFICATION_2} combined with Proposition \ref{DUALITY_EXTENDED} provides (non-canonical)
complex isomorphisms
\begin{enumerate}
\item[{\rm(1)}]
 $\mathcal D_{\phi}(t)\iso \mathcal C_{\phi^*}(s-t)$ and
\item[{\rm(2)}] $ \mathcal C_{\psi}(t)\iso \mathcal
D_{\psi^*}(r-t)$
\end{enumerate}
where as above $s=n-l,\ r=n-m$. Next we consider the diagram
$$
\CD \mathcal D_{\phi}(t+m+1)\otimes
F^*@>\partial_{\psi}>>\mathcal D_{\phi}(t+m)@>>>\mathcal M(t)@>>>0\\
@VVV @VVV\\
\mathcal C_{\phi^*}(\rho-t-1)\otimes F^*@>d_{\psi^*}>>\mathcal
C_{\phi^*}(\rho-t)@>>>\mathcal C_{\bar{\lambda}}(\rho-t)@>>>0.
\endCD
$$
The isomorphisms (1) assure that the vertical arrows are
isomorphisms, while (2) provides the commutativity of the diagram.
The desired isomorphism is induced.
\end{proof}

From the Koszul bicomplex in the previous section we extract the
sequence of complexes
$$
\CD \mathcal D_{\phi}(t+m+1)\otimes
F^*@>\partial^{-1}_{\psi}>>\mathcal
D_{\phi}(t+m)@>\pm\nu_\psi>>\mathcal
D_{\phi}(t)@>\partial^{0}_{\psi}>>\mathcal D_{\phi}(t-1)\otimes F.
\endCD
$$
Since $\mathcal M(t)=\Coker (B^{t,-1}\to B^{t,0})$ and $\mathcal
N(t)=\Ker (C^{t,0}\to C^{t,1})$, we obtain an induced complex
morphism
$$
\CD \mathcal M(t)@>\nu>>\mathcal N(t).
\endCD
$$
This, composed with the inverse of the isomorphism of Proposition
\ref{N identification2}, yields a morphism
$$
\CD \mathcal C_{\bar{\lambda}}(\rho-t)@>\mu>>\mathcal N(t).
\endCD
$$
Using Theorem \ref{E-N HOMOLOGY} we easily deduce the following
properties of $\mu$.

\begin{proposition}\label{GRADESENSIBILITY} Let $r,l,g$ be as above.
\begin{enumerate}
\item[\rm{(1)}]Suppose that $t+l\le 0$ or $r<t+l$. Then the
$\mu_i$ are isomorphisms for $i>r+1-g$, and $\mu_{r+1-g}$ is
injective. \item[\rm{(2)}]Suppose that $l\le t+l \le r$.
\begin{enumerate}
\item[\rm{(i)}]If $r+1-g\le t$, then the $\mu_i$ are isomorphisms
for $i>r+1-g$, and $\mu_{r+1-g}$ is injective. \item[\rm{(ii)}]If
$t+l\le r+1-g$, then the $\mu_i$ are isomorphisms for $i>r+2-g-l$,
and $\mu_{r+2-g-l}$ is injective. \item[\rm{(iii)}]If
$t<r+1-g<t+l$, then the $\mu_i$ are isomorphisms for $i>t$.
\end{enumerate}
\item[\rm{(3)}]Suppose that $0<t+l<l$. Then the $\mu_i$ are
isomorphisms for $i>\min(0,r+1-g-l-t)$ and, if $r+1-g-l-t\ge 0$,
then $\mu_{r+1-g-l-t}$ is injective.
\end{enumerate}
\end{proposition}

\begin{theorem}\label{MaximalCaseHomology} Suppose that $g\ge r+1$. With notation as above set $S_0(C)=R/I_\chi$.
Equip $\mathcal C_{\bar{\lambda}}(t)$ with the graduation induced
by the complex isomorphism $\mathcal M(\rho -t)\to \mathcal
C_{\bar{\lambda}}(t)$ of Proposition $\ref{N identification2}$.
Then for the homology $\tilde H^{^.}$ of $\mathcal
C_{\bar{\lambda}}(t)$ the following holds:
\begin{enumerate}
\item[\rm{(a)}]in case $t\le \frac{\rho}2$,
$$
\tilde H^{i}=
\begin{cases}
D_{\rho-t-\frac{i-1}2}(\mathcal{H}^*)\otimes S_{\frac{i-1}2}(C) &\quad\text{if}\ \ 0\le i< h,\; i\not\equiv 0\ (2),\\
\phantom{D_{\rho-t-\frac{i-1}2}(\mathcal{H}^*)}0 &\quad\text{if}\
\ 0\le i< h,\; i\equiv 0\ (2);
\end{cases}
$$
\item[\rm{(b)}]in case $\frac{\rho}2<t\le\rho$, \small
$$
\tilde H^{i}=
\begin{cases}
D_{\rho-t-\frac{i-1}2}(\mathcal{H}^*)\otimes S_{\frac{i-1}2}(C) &\quad\text{if}\ \ 0\le i< \min(h,2(\rho-t+1)),\; i\not\equiv 0\ (2),\\
S_{\frac{i-l}2-\rho+t-1}(\mathcal{H})\otimes S_{\frac{i+l}2-1}(C) &\quad\text{if}\ \ 2(\rho-t+1)+l\le i< h,\; i-l\equiv 0\ (2),\\
\phantom{D_{\rho-t-\frac{i-1}2}(\mathcal{H}^*)}0
&\quad\text{otherwise if}\ \ 0\le i< h;
\end{cases}
$$
\normalsize \item[\rm{(c)}]in case $\rho<t<r$,
$$
\tilde H^{i}=
\begin{cases}
S_{\frac{i-r+t-1}2}(\mathcal{H})\otimes S_{\frac{i+r-t-1}2}(C) &\quad\text{if}\ \ r-t+1\le i< h,\; i+r-t\not\equiv 0\ (2),\\
\phantom{S_{\frac{i-r+t-1}2}(\mathcal{H})}0 &\quad\text{otherwise
if}\ \ 0\le i< h;
\end{cases}
$$
\normalsize
 \item[\rm{(d)}]in case $r\le t$,
$$
\tilde H^{i}=
\begin{cases}
S_{\frac{i-1}2+t-r}(\mathcal{H})\otimes S_{\frac{i-1}2}(C) &\quad\text{if}\ \ 0\le i< h,\; i\not\equiv 0\ (2),\\
\phantom{S_{\frac{i-r+t-1}2}(\mathcal{H})}0 &\quad\text{if}\ \
0\le i< h,\; i\equiv 0\ (2).
\end{cases}
$$
\end{enumerate}
\end{theorem}

\begin{proof} If $h=0$, then there is nothing to prove. If $h\ge 1$ then $\rho\ge 0$ and
$r\ge 1$.

We use the complex morphism $\mu:\mathcal
C_{\bar{\lambda}}(t)\to\mathcal{N}_.(\rho-t)$ from above which
induces the following commutative diagram:
$$
\begin{CD}
C^{-1}_{\bar{\lambda}}@>\partial^{-1}_{\bar{\lambda}}>>C^{0}_{\bar{\lambda}}@>\partial^{0}_{\bar{\lambda}}>>C^{1}_{\bar{\lambda}}@>>> C^{2}_{\bar{\lambda}}\\
&& @V\mu_0 VV  @V\mu_{1} VV @V\mu_{2}VV \\
0 @>>> N^0@>\bar d^0_\phi>>N^1 @>>> N^2\\
&& @VVV  @VVV @VVV\\
0@>>> \Coker\mu_0@>>>  0@>>>  0
\end{CD}
$$
Since $C^{-1}_{\bar{\lambda}}$ is a torsion module and $\mu_0$ is
injective, we have $\partial^{-1}_{\bar{\lambda}}=0$. If $h\ge 1$,
then $d^0_\phi$ is injective. This implies that
$\partial^{0}_{\bar{\lambda}}$ is injective, so $\tilde H^{0}=0$.
If $h\ge 2$, Theorem \ref{FUNDAMENTAL} (or Theorem
\ref{FUNDAMENTALnegativ}) says that the row homology at $N^0$ and
at $N^1$ is 0, so
$$
\tilde H^{1}=\Coker\mu_0=
\begin{cases}
D_{\rho-t}(\mathcal{H}^*)\otimes R/I_\chi &\quad\text{if}\ \ t\le\rho,\; \\
\phantom{D_{\rho-t}(\mathcal{H}^*)}0 &\quad\text{if}\ \ \rho< i< r,\;\\
S_{t-r}(\mathcal{H})\otimes R/I_\chi &\quad\text{if}\ \ r\le t. \\
\end{cases}
$$

If $h\ge 3$, then $\tilde H^{2}$ equals the row homology at $N^2$.
The remaining statements follow easily from Theorem
\ref{FUNDAMENTAL_2} (or Theorem \ref{FUNDAMENTALnegativ}).
\end{proof}

\begin{remark}The case $t=0$ in Theorem \ref{MaximalCaseHomology}
covers Proposition 5.1 in [MNP].
\end{remark}

\begin{corollary}\label{Case l Big}Suppose that $h=|\, \rho\, |+1$.
\begin{enumerate}
\item[\rm{(1)}]If $l\ge \frac{r-1}2$, then $\big(\mathcal
C_{\bar{\lambda}}(0)\big)^*$ has non-vanishing homology only in
degree $\rho+1$, and if $r>1$, then
$$\bar H^{0}(\rho+1)=S_{\rho}(\Coker\chi^*).$$
\item[\rm{(2)}]If $l\ge \frac{r+1}2$, then the homology of
$C^{^.}_{\bar{\lambda}}(\rho+1)$ vanishes in positive degrees
except for $\rho+1$, and
$$\tilde H^{\rho+1}(\rho+1)=S_{\rho+1}(\Coker\chi^*).$$
\end{enumerate}
\end{corollary}

\begin{proof}(a) Since $h\ge 1$, we have $\rho\ge 0$. In positive degrees, the homology of
$\big(\mathcal C_{\bar{\lambda}}(0)\big)^{*}$ is almost the same
as the homology of $\mathcal C_{\bar{\lambda}}(\rho)$, with the
only exception in grade $1$ (where the homology of $\big(\mathcal
C_{\bar{\lambda}}(0)\big)^{*}$ is 0). If we require that $l\ge
\frac{r-1}2$, Proposition \ref{GRADESENSIBILITY},(2) provides the
result.
\bigskip

(b) If $r=1$, then the claim is clear. If $r>1$, then $l>1$, and
the result follows directly from Proposition
\ref{GRADESENSIBILITY},(3)
\end{proof}

\begin{remark}Corollary \ref{Case l Big},(1) was inspired by Lemma 5.5 in
[MPN]. Corollary \ref{Case l Big},(2) can be seen as an extension
of Theorem \ref{E-N HOMOLOGY},(b), since $l$ being big enough, one
can easily deduce similar results.
\end{remark}

\begin{remark}In the last part of this section we required for $\grade I_\chi$ to have the
greatest possible value. If we do the same for $\grade I_\lambda$,
we can obtain information about the homology of $\mathcal
C_{\bar{\lambda}}(t)$, by studying the upper half of the Koszul
bicomplex $\mathcal{K}_{.,.}(t)$.
\end{remark}


\begin{thebibliography}{15.}
\addcontentsline{toc}{section}{Bibliography}

\bibitem[BO1]{BO1} N.Bourbaki. {\em Algebra I, Chapters 1-3},
Springer 1989.

\bibitem[BO2]{BO2} N.Bourbaki. {\em Alg\`{e}bre I, Chapitre 10, Alg\`{e}bre homologique}, Masson 1980.

\bibitem[BE1]{BE1} D. Buchsbaum and D. Eisenbud. {\em Some Structure Theorems for Finite
Free Resolutions}. Advances in Math. {\bf 12} (1974), 84--139.

\bibitem[BE2]{BE2} D. Buchsbaum and D. Eisenbud. {\em Generic Free Resolutions and a Family of
Generically Perfect Ideals}. Advances in Math. {\bf 18} (1975),
245--301.

\bibitem[BE3]{BE3} D. Buchsbaum and D. Eisenbud. {\em Remarks on Ideals and Resolutions}. Symposia Math. {\bf 11} (1973), 193--204.

\bibitem[BH]{BH} W.Bruns and J. Herzog. {\em Cohen-Macaulay Rings}. Cambridge Uni. Press 1996.

\bibitem[BV1]{BV1} W. Bruns and U. Vetter. {\em Length formulas for the Local
Cohomology of Exterior Powers}. Math Z. {\bf 191} (1986),
145--158.

\bibitem[BV2]{BV2} W. Bruns and U. Vetter. {\em Determinantal rings}.
Lect.~Notes Math. {\bf 1327}, Springer 1988.

\bibitem[BV3]{BV3} W. Bruns and U. Vetter. {\em A Remark on Koszul Complexes}.
Beitr. Algebra Geom. {\bf 39} (1998), 249--254.

\bibitem[BV4]{BV4} W. Bruns and U. Vetter. {\em The Koszul Complex in Projective Dimension One}.
Geometric and Combinatorial Aspects of Commutative Algebra, J.
Herzog and G. Restuccia eds, Marcel Dekker 2001, 89--98.

\bibitem[E]{E} D. Eisenbud. {\em Commutative Algebra with a View Toward Algebraic Geometry}.
Grad. Text in Math. {\bf 150}, Springer 1995.

\bibitem[HM]{HM} J. Herzog and A. Martsinkovsky. {\em Glueing Cohen-Macaulay
modules with applications to quasihomogeneous complete
intersections with isolated singularities}. Comment. Math. Helv.
{\bf 68} (1993), 365--384.

\bibitem[I]{I} B. Ichim. {\em Generalized Koszul Complexes}. Thesis,
Universität Oldenburg (Germany), 2004.

\bibitem[MNP]{MNP} J.C. Migliore, U. Nagel, and C. Peterson. {\em Buchsbaum-Rim
sheaves and their multiple sections}. J. Algebra {\bf 219} (1999),
378--420.

\bibitem[N]{N} D.G. Northcott. {\em Finite free resolutions}.
Cambridge Uni. Press {\bf 71}, 1976.


\end{thebibliography}
\end{document}